\documentclass[reqno]{amsart}
\usepackage{amsfonts}

\usepackage{amssymb}
\usepackage{xcolor}
\usepackage[colorlinks=true]{hyperref}

\hypersetup{urlcolor=blue, citecolor=red}
\newtheorem{Lemma}{Lemma}[section]
\newtheorem{Theorem}[Lemma]{Theorem}

\newtheorem{Corollary}[Lemma]{Corollary}

\begin{document}
\title[An inverse kinematic problem with internal sources]{An inverse
kinematic problem with internal sources}
\author[L. Pestov]{Leonid Pestov}
\address{Immanuel Kant Baltic Federal University, Russia}
\email{LPestov@kantiana.ru}
\author[G. Uhlmann]{Gunther Uhlmann}
\address{Department of Mathematics, University of Washington, USA}
\email{gunther@math.washington.edu}
\author[H. Zhou]{Hanming Zhou}
\address{Department of Mathematics, University of Washington, USA}
\email{hzhou@math.washington.edu}

\begin{abstract}
Given a bounded domain $M$ in $\mathbb{R}^n$ with a conformally Euclidean
metric $g=\rho\,dx^2$, in this paper we consider the inverse problem of
recovering a semigeodesic neighborhood of a domain $\Gamma\subset \partial M$
and the conformal factor $\rho$ in the neighborhood from the travel time
data (defined below) and the Cartesian coordinates of $\Gamma$. We develop
an explicit reconstruction procedure for this problem. The key ingredient is
the relation between the reconstruction and a Cauchy problem of the
conformal Killing equation.
\end{abstract}

\maketitle


\section{Introduction}

Let $(M,g)$ be a bounded domain in $\mathbb{R}^n$, $n\geq 2$ with smooth
boundary $\partial M$. We assume $g$ is conformal to the Euclidean metric,
i.e. $g=\rho\,dx^2$, where $\rho$ is a positive smooth function on $M$ and $%
dx^2=(dx^1)^2+\cdots+(dx^n)^2$ is the Euclidean metric. Let $\Gamma$ be a
domain in $\partial M$ (in particular $\Gamma$ could be $\partial M$), from
each $x^{\prime }\in \Gamma$, there is a unique geodesic $\gamma_{x^{\prime
}}(t)$ with $\gamma_{x^{\prime }}(0)=x^{\prime },\, \dot{\gamma}_{x^{\prime
}}(0)=\nu(x^{\prime })$, where $\nu(x^{\prime })$ is the inward unit normal
vector to $\partial M$ at $x^{\prime }$ w.r.t. the metric $g$. Moreover,
since $\gamma$ is a geodesic of unit speed w.r.t. the metric $g$, we have $%
\rho(\gamma)|\dot{\gamma}|^2=1$, where $|\cdot|$ is the Euclidean norm.

There is a positive smooth function $T(x^{\prime })$ on $\Gamma$ such that
for each $x^{\prime }\in\Gamma$, the geodesic $\gamma_{x^{\prime }}$, which
is orthogonal to $\partial M$ at $x^{\prime }$, is defined on the interval $%
[0,T(x^{\prime })]$. Let 
\begin{equation*}
D:=\{(x^{\prime },t): x^{\prime }\in\Gamma,\, 0\leq t< T(x^{\prime })\},
\end{equation*}
we consider the map $\gamma: D\to \gamma(D),\, \gamma(x^{\prime
},t):=\gamma_{x^{\prime }}(t)=x$. Generally such a map is not a
diffeomorphism, for example given the Euclidean disk with radius $r$, let $%
\Gamma$ be a domain of the boundary, then $\gamma$ is not a diffeomorphism
if $T(x^{\prime })> r$. Thus we modify $T(x^{\prime })$, i.e. $D$, so that $%
\gamma: D\to\gamma(D)$ is a diffeomorphism. Under this assumption, $D$
actually provides a semigeodesic coordinate system (or boundary normal
coordinates) for $\gamma(D)$.

Now given a point $x\in \gamma (D)$, if $x^{\prime }\in \Gamma $ such that $%
x=\gamma (x^{\prime },t)$ for some $t\in \lbrack 0,T(x^{\prime }))$, let $%
U(x^{\prime })\subset \Gamma $ be a neighborhood of $x^{\prime }$. Moreover,
we fix $U(x^{\prime })$ for all $x\in \gamma _{x^{\prime }}([0,T(x^{\prime
})))$, and $U(x^{\prime })$ can be arbitrarily small. Notice that the choices of $U(x^{\prime})$ somehow depend on a priori knowledge of $x$, i.e. we need to know the geodesic projection $x^{\prime}(x)\in \Gamma$ of any $x\in \gamma(D)$. We define the 
\textit{travel time data} w.r.t. $D$ by 
\begin{equation*}
\Omega (D):=\{(\tau (x,x^{\prime \prime }), x^{\prime}(x)) : x\in \gamma (D),\,x^{\prime \prime }\in
U(x^{\prime })\},
\end{equation*}
where $\tau(x,x^{\prime \prime }):=dist_{g}(x,x^{\prime \prime})$. In this paper, we consider the problem of recovering the neighborhood $%
\gamma (D)$ and the conformal factor $\rho $ in $\gamma (D)$ from the travel
time data $\Omega (D)$. To solve the problem, we need some extra inverse
data, namely we assume the Cartesian coordinates of $\Gamma $ is known. This
is a reasonable assumption, since any rigid transformation of the domain $M$
does not change the travel time data. We call the problem the \textit{%
inverse kinematic problem with internal sources}. It's worth pointing out
that we do not put any assumption on the convexity of the boundary (or 
$\Gamma$). Uniqueness for this inverse problem was proved by
Yu. E. Anikonov \cite{An}. In this paper we give a reconstructon
procedure, based on conformal Killing vector fields.

\begin{Theorem}
\label{main} Let $M$ be a bounded domain in $\mathbb{R}^n$, $n\geq 2$ and $%
g=\rho\,dx^2$ be a conformally Euclidean metric on $M$. Let $\Gamma$ be a
domain in $\partial M$, then there exists a semigeodesic coordinate system $D
$ such that $\gamma:D\to \gamma(D)\subset\mathbb{R}^n$ is a diffeomorphism
and $\Gamma=\gamma(\{t=0\})$. We develop a reconstruction procedure of the
diffeomorphism $\gamma$ and the conformal factor $\rho$ in $\gamma(D)$ from
the travel time data $\Omega(D)$ and the Cartesian coordinates of $\Gamma$.
\end{Theorem}

Generally one can not expect to reconstruct $\rho$ on the whole manifold, as
the necessary assumption that $\gamma$ is a diffeomorphism. However, if $%
M\backslash \gamma(D)$ has empty interior, we reconstruct the domain $M$ and
the conformal factor $\rho$ globally by taking limit. The example mentioned
above satisfies the assumption, and it is easy to see that $%
M\backslash\gamma(D)$ is the center of the disk if $\Gamma=\partial M,\,
T(x^{\prime })\equiv r$.

\begin{Corollary}
Let $M$ be a bounded domain in $\mathbb{R}^n$, $n\geq 2$ and $g=\rho\,dx^2$
be a conformally Euclidean metric on $M$. Assume that there exists a
semigeodesic coordinate system $\mathcal{D}=\{(x^{\prime },t):x^{\prime
}\in\partial M, t\in[0,T(x^{\prime }))\}$ such that $\gamma:\mathcal{D}\to
\gamma(\mathcal{D})\subset\mathbb{R}^n$ is a diffeomorphism with $%
\gamma(\{t=0\})=\partial M$, and $M\backslash\gamma(\mathcal{D})$ has empty
interior. Then there is a reconstruction procedure of the domain $M$ and the
conformal factor $\rho$ from the travel time data $\Omega(\mathcal{D})$ and
the Cartesian coordinates of $\partial M$.
\end{Corollary}

The inverse kinematic problem arose in geophysics in an attempt to determine
the substructure of the Earth by measuring at the surface of the Earth the
travel times of seismic waves. In application the conformal factor $\rho$
corresponds to $1/c^2(x)$, where $c(x)$ is the sound speed (index of
refraction). It goes back to \cite{H1905, WZ1907} who considered the case of
a radial metric conformal to the Euclidean metric. The case considered above
corresponds to an isotropic media, however it has been realized later that
the inner core of the Earth might exhibit anisotropic behavior \cite{Cr92}.
The geometric version of the problem is related to the boundary rigidity
problem and its linearization, namely the geodesic ray transform, see \cite%
{SU08} for a recent survey.

However, in current paper we also take use of the internal data (data from
internal sources), so that we can reconstruct the geometry explicitly. In 
\cite{Ku97, KKL01} an isometric copy of a compact Riemannian manifold was recovered
from the set of boundary distance functions $\{r_{x}(y):=dist(x,y):x\in
M,y\in \partial M\}$, which is also internal data. Such internal data is
also related to the broken geodesic flow which consists of two geodesic
segments sharing a common end point inside the manifold, see e.g. \cite%
{KKL01, KLU10}. A related reconstruction problem with different assumptions was considered in \cite{DHILU} by reducing the travel time data to measurements of the shape operators of the wave fronts of waves diffracted from interior points. Different from our method, their approach treated the case of two dimensions and the case of three and higher dimensions separately.

Notice that the statement of Theorem \ref{main}
actually shows that this is a local problem, i.e. for a point $x=\gamma
(x^{\prime },t)$, we only need the travel time data from $x$ to an
arbitrarily small neighborhood $U(x^{\prime })\subset \Gamma $. If the
function $T(x^{\prime })$ is also uniformly small, the problem can be
formulated just near one boundary point. Similar to the local version of the
problem, the local boundary rigidity problem and local geodesic ray
transform were considered in \cite{UV12, SUV13}, and a generalization to local ray transforms along arbitrary smooth curves was studied in \cite{Zh13}. 

As mentioned above, our arguments give a reconstruction procedure for the
diffeomorphism $\gamma $ and the conformal factor $\rho $. The
reconstruction procedure consists of two steps: step 1 (Section 2) devotes
to the recovery of a semigeodesic (isometric) copy of the metric $g$ and the
boundary restriction of the conformal factor from the inverse data; in step
2 (Section 3), we reconstruct $\gamma $ and $\rho $ by studying the relation
between $\gamma $ and conformal Killing vectors on the semigeodesic copy $D$.

Notice also that the inverse dynamical problem for the wave equations
(with boundary data) may be reduced by the boundary control method to the
inverse kinematic problem with internal sources and then to the Yamabe
problem \cite{BG01}. It gives in our case the Cauchy problem for
the Laplace operator. We use another approach (using conformal Killing
vector fields) that gives stability for dimensions $n>2$ \cite
{R94}. 

$\newline
$\noindent \textbf{Acknowledgements.} L. Pestov was partly supported by
grant 12-01-00260-a, RFBR; G. Uhlmann was partly supported by NSF and
a Simons Fellowship; H. Zhou was partly supported by NSF.


\section{Recovery of the semigeodesic copy of the metric}

Given $x\in\gamma(D)$, there is a unique geodesic $\gamma_{x^{\prime }}$
(normal to $\partial M$) and $0\leq t<T(x^{\prime })$ such that $%
x=\gamma(x^{\prime },t)$. Here $x^{\prime }(x)$ is the geodesic projection
of $x$ on $\Gamma$, so $t(x)=\tau(x,x^{\prime }(x))$ is the distance from $x$
to the boundary. 
Thus the travel time data $\Omega(D)$ uniquely determines the semigeodesic coordinates of points in $\gamma(D)$.

Let the pair $(x^{\prime }(x),t(x))$ be the semigeodesic coordinate of the
point $x$. Thus for any point $y=(x^{\prime },t)\in D$, the function $%
\tau(\gamma(x^{\prime },t), x^{\prime \prime })$ is known. Define 
\begin{equation*}
\lambda((x^{\prime },t),x^{\prime \prime }):=\tau(\gamma(x^{\prime
},t),x^{\prime \prime }),
\end{equation*}
then $\lambda(y,x^{\prime \prime }),\, y=(x^{\prime },t)$ is the distance
between points $y\in D$ and $x^{\prime \prime }\in U(x^{\prime })$ in the
metric $\tilde{g}:=\gamma^*(g)$. Note that we identify $\Gamma$ with $%
\Gamma\times\{0\}$. We call $\tilde{g}$ the \textit{semigeodesic copy} of
the metric $g$.

We first recover the metric $\tilde{g}$. Notice that $\gamma$, now as an
isometry, sends geodesics to geodesics, this implies 
\begin{equation*}
\tilde{g}_{kn}(x^{\prime },0)=\delta_{kn},\, x^{\prime }\in\Gamma,\, 1\leq
k\leq n,
\end{equation*}
where $\delta_{ij}$ is the Kronecker delta. In local coordinates, $%
y=(y^1,\cdots, y^n)$, where $y^n=t$, one has the following eikonal
equation 
\begin{equation}  \label{eikonal}
1=|\nabla\tau(x,x^{\prime \prime })|_g^2=|\tilde{\nabla}\lambda(y,x^{\prime
\prime })|_{\tilde{g}}^2=\tilde{g}^{ij}(y)\frac{\partial\lambda(y,x^{\prime
\prime })}{\partial y^i}\frac{\partial\lambda(y,x^{\prime \prime })}{%
\partial y^j},\, x^{\prime \prime }\in U(x^{\prime }).
\end{equation}
Note that if $y^{\prime }$ is close enough to $y$, and $x^{\prime \prime
}\in U(x^{\prime }(y))$ sufficiently close to $x^{\prime }(y)$, then $%
x^{\prime \prime }\in U(x^{\prime }(y^{\prime }))$ too. (\ref{eikonal})
gives a family of linear algebraic equations w.r.t. the contravariant
components $\tilde{g}^{ij}(y)$ of the metric $\tilde{g}$. For $y\in D$ with $%
y^n>0$, it is known that $\tilde{g}^{ij}(y)$ can be recovered by the
knowledge of $\tilde{g}^{ij}(y)\partial_{y^i}\lambda(y,x^{\prime \prime
}_k)\partial_{y^j}\lambda(y,x^{\prime \prime }_k)$ for $N=n(n+1)/2$
``generic" points $x^{\prime \prime }_k\in U(x^{\prime }(y)),\,
k=1,2,\cdots, N$, see e.g. \cite{Sh94, SU98}. Such $N$ generic points always
exist in a neighborhood of $x^{\prime }(y)$ in $\Gamma$. Thus $\tilde{g}%
^{ij}(y)$ is determined for $y\in D,\, y^n>0$. By the smoothness of $\tilde{g%
}$, it also recovers $\tilde{g}^{\alpha\beta}(x^{\prime },0),\,1\leq
\alpha,\beta\leq n-1,\, x^{\prime }\in\Gamma$. So the metric $\tilde{g}$ is
uniquely determined.

Notice that the boundary restriction of the isometry $\gamma$ is given, i.e. 
$x^{\prime }=x^{\prime}(y^{1},\cdots,y^{n-1},0)=(x^{\prime 1},\cdots, x^{\prime n})$ as a point in $\mathbb{R}^{n}$ is known, we recover the conformal factor $\rho$ on $%
\Gamma$. Since $\tilde{g}=\gamma^*(g)$ is the pullback, we have 
\begin{equation}  \label{pullback}
\tilde{g}_{\alpha\beta}(x^{\prime },0)=\frac{\partial\gamma^k}{\partial
y^{\alpha}}\frac{\partial\gamma^l}{\partial y^{\beta}}\rho(x^{\prime
})\delta_{kl}=\rho(x^{\prime })\sum_{k=1}^n\frac{\partial x^{\prime k}}{%
\partial y^{\alpha}}\frac{\partial x^{\prime k}}{\partial y^{\beta}},\,
\alpha, \beta=1,\cdots,n-1.
\end{equation}
Thus equation (\ref{pullback}) and the knowledge of $\tilde{g}|_{t=0}$
together determine $\rho|_{\Gamma}$.


\section{Recovery of the isometry and the conformal factor}

To recover the conformal factor $\rho$, we need to solve the \textit{%
pullback problem}, i.e. to find the map $\gamma$. To this end, we need some
knowledge of conformal Killing vector fields. Recall that a vector field $u$
is called a conformal Killing vector field if it satisfies the conformal
Killing equation ( in covariant form ) 
\begin{equation*}
Ku:=\sigma\nabla u-g\delta u/n=0,
\end{equation*}
where $\sigma\nabla$ is the symmetric part of the covariant derivative $%
\nabla$, $\delta$ is the divergence. They are exactly those vector fields
whose flows preserve the conformal structures of the manifolds. In local
coordinates, the conformal Killing equation has the form ( for covariant
components ) 
\begin{equation*}
\frac{1}{2}(\frac{\partial u_i}{\partial x^j}+\frac{\partial u_j}{\partial
x^i}-2\Gamma_{ij}^ku_k)-\frac{1}{n}g_{ij}(g^{kl}\frac{\partial u_l}{\partial
x^k}-g^{kl}u_m\Gamma_{kl}^m)=0.
\end{equation*}
Thus the conformal Killing equation is equivalent to a system of first order
partial differential equations.

However, not all of the metrics admit conformal Killing vector fields,
actually for $n\geq 3$, a ``generic'' metric does not possess any
non-trivial conformal Killing vector fileds, see e.g. \cite{BCS05, LS12}. On
the other hand, note that the metric $g$ is conformal to the Euclidean
metric, thus they share the same set of conformal Killing vector fields.
When $n=2$, in the Cartesian coordinates $(x^1,x^2)$ all conformal Killing
vector fields of the Euclidean metric have the form $u=(u^1, u^2)$, where $%
u^1$ and $u^2$ are conjugate harmonic functions. In the case $n>2$, the
contravariant components of $u$ in the Cartesian coordinates $(x^1,\cdots,
x^n)$ are given by 
\begin{equation*}
u^i(x)=a_0x^i+(Ax)^i-b^i|x|^2+2x^i(b,x)+c^i,
\end{equation*}
where $a_0$ is a real constant, $A$ is a $n\times n$ skew-symmetric constant
matrix, $b$ and $c$ are vectors in $\mathbb{R}^n$.

Now we are in a position to recover the map $\gamma$ and the conformal
factor $\rho$. Let $e_{(j)}=\frac{\partial}{\partial x^j},\, j=1,\dots,n$ be
the standard basis vectors in $\mathbb{R}^n$. It is easy to see that they
are conformal Killing vector fields in Euclidean metric, thus also conformal
Killing vector fields for $g$. Then $u_{(j)}=\gamma^* e_{(j)},\, j=1,\dots,n$
are conformal Killing vector fileds for the metric $\tilde{g}=\gamma^* g$.
This implies $u_{(j)},\, j=1,\dots,n$ satisfy the conformal Killing equation 
\begin{equation*}
Ku_{(j)}(y)=0,\, y\in D,\, j=1,\dots,n.
\end{equation*}
It is known that $u_{(j)}$ is uniquely determined by the Cauchy data $%
\{u_{(j)}(x^{\prime },0): x^{\prime }\in\Gamma\}$, see e.g. \cite{L97, DS11}%
. Thus we calculate the Cauchy data of $u_{(j)}$ first.

Since we have already recovered the semigeodesic copy $\tilde{g}$, we denote
the dual vector of $u_{(j)}$ by $u^{(j)}$. In the mean time, we denote the
dual vector of $e_{(j)}$ under the metric $g=\rho\,dx^2$ by $w^{(j)}$, then $%
w^{(j)}=\rho\,dx^j$. In local coordinates, the equality $u_{(j)}=\gamma^*
e_{(j)}$ means ( for covariant components ) 
\begin{equation*}
u^{(j)}_i(y)=w^{(j)}_k\frac{\partial x^k}{\partial y^i}=\rho\frac{\partial
\gamma^j(y)}{\partial y^i}.
\end{equation*}
This observation is crucial in our reconstruction procedure, it relates the
isometry $\gamma$ to the conformal Killing vector fields on the semigeodesic
copy $D$. Thus at $y^n=t=0$, 
\begin{equation*}
u^{(j)}_{\alpha}(x^{\prime },0)=\rho(x^{\prime })\frac{\partial\gamma^j}{%
\partial y^{\alpha}}(x^{\prime },0)=\rho(x^{\prime })\frac{\partial
x^{\prime j}}{\partial y^{\alpha}},\quad \alpha=1,\cdots, n-1.
\end{equation*}
Since $\rho(x^{\prime })$ and $\gamma(x^{\prime },0)$ are known for $%
x^{\prime }\in\Gamma$, $u^{(j)}_{\alpha}(x^{\prime },0)$ are determined. To
determine the value of $u^{(j)}_n$ at $t=0$, notice that $\dot{\gamma}%
_{x^{\prime }}(0)=\nu(x^{\prime })$, so 
\begin{equation*}
u^{(j)}_n(x^{\prime },0)=\rho(x^{\prime })\frac{\partial \gamma^j}{\partial t%
}(x^{\prime },0)=\rho(x^{\prime})\nu^{j}(x^{\prime }).
\end{equation*}
However, $\nu$ is a unit vector w.r.t. metric $g=\rho\,dx^2$, if we denote
the inward unit normal vector on $\partial M$ w.r.t. the Euclidean metric by 
$\nu_0$, then $\nu=\frac{1}{\sqrt{\rho}}\nu_0$, i.e. 
\begin{equation*}
u^{(j)}_n(x^{\prime },0)=\sqrt{\rho(x^{\prime })}\nu_0^j(x^{\prime }).
\end{equation*}
Fortunately, the Cartesian coordinates of the hypersurface $\Gamma$ are
given, thus $\nu_0$ as the normal vector to $\Gamma$ is known. Together with
the knowledge of $\rho|_{\Gamma}$, we recover $u^{(j)}_n|_{t=0}$ too.

From the Cauchy data, we uniquely recover the conformal Killing vector
fields $u_{(j)}=\gamma^*e_{(j)}$, equivalently the dual vector $u^{(j)}$.
Now by defining $v=(v^1,\cdots, v^n)$, 
\begin{equation*}
v^j(x^{\prime},t):=u^{(j)}_n(x^{\prime },t)=\rho(\gamma(x^{\prime },t))\frac{%
\partial \gamma^j(x^{\prime },t)}{\partial t},
\end{equation*}
we have $v=\rho(\gamma)\,\dot{\gamma}$. In the mean time, notice that 
\begin{equation}  \label{rho}
|v|^2=\rho^2(\gamma)|\dot{\gamma}|^2=\rho(\gamma)\quad (\mbox{since}\, |\dot{%
\gamma}|_g=1),
\end{equation}
we obtain 
\begin{equation*}
\dot{\gamma}=\frac{v}{|v|^2}.
\end{equation*}
This implies 
\begin{equation*}
\gamma(x^{\prime },t)=\int_0^t\dot{\gamma}_{x^{\prime }}(t)\,dt+x^{\prime
}=\int_0^t \frac{v}{|v|^2}(x^{\prime },t)\, dt+x^{\prime },
\end{equation*}
i.e. we recover the geodesics $\gamma_{x^{\prime }}(t)$, therefore the
diffeomorphism $\gamma: D\to \gamma(D)$, namely the range $\gamma(D)$.
Moreover, by (\ref{rho}) 
\begin{equation*}
\rho(\gamma(x^{\prime },t))=|v(x^{\prime},t)|^{2},
\end{equation*}
the conformal factor $\rho|_{\gamma(D)}$ is determined.


\end{document}